\documentclass{amsart}

\usepackage{amssymb}
\usepackage{eepic}
\usepackage{color}

\allowdisplaybreaks

\numberwithin{equation}{section}

\newtheorem{theorem}{Theorem}[section]
\newtheorem{lemma}[theorem]{Lemma}

\newtheorem{remark}[theorem]{Remark}

\newtheorem{TheoA}{Theorem A}
\newtheorem{TheoB}{Theorem B}

\newcommand{\N}{\mathbb{N}}
\newcommand{\T}{\mathbb{T}}
\newcommand{\Z}{\mathbb{Z}}
\newcommand{\R}{\mathbb{R}}
\newcommand{\C}{\mathbb{C}}

\def\1{\mathbf{1}}

\newcommand{\dem}{\noindent {\bf Proof. }}
\newcommand{\ske}{\noindent {\bf Sketch of the proof. }}
\newcommand{\demB}{\noindent {\bf Proof of Theorem B. }}
\newcommand{\fin}{\hspace*{\fill} $\square$ \vskip0.2cm}

\def\mean{- \hskip-10.5pt \int}

\addtocounter{tocdepth}{-1}

\begin{document}

\addtolength{\parskip}{+1ex}

\title[On the growth of vector-valued Fourier series]{On the growth of \\ vector-valued Fourier series}

\author{Javier Parcet, Fernando Soria and Quanhua Xu}

\maketitle



\begin{abstract}
Let $f: \mathbb{T} \to \mathrm{X}$ satisfy $$\int_\mathbb{T} \|f(x)\|_\mathrm{X} \big( \log^+ \|f(x)\|_\mathrm{X} \big)^{1+\delta} \, dx \, < \, \infty,$$ where $\mathrm{X}$ is a UMD Banach space and $\delta > 0$. Then we prove that $$\Big\| \sum_{|k| \le n} \widehat{f}(k) e^{2\pi i k x} \Big\|_\mathrm{X} \, = \, o \big( \log \log n \big) \quad \mbox{for} \quad \mathrm{ae} \hskip1pt - \hskip1pt x \in \T.$$ In other words, the \lq little Carleson theorem\rq${}$ holds for UMD-valued functions.  
\end{abstract}

\vskip1cm

\section*{Introduction}

In 1966, L. Carleson proved \cite{C} that the Fourier series of any square integrable function $f: \T \to \C$ converges almost everywhere to $f$. This result is a corner stone in the harmonic analysis of the 20th century. Over the years, Carleson's theorem was reproved by Fefferman, Lacey/Thiele or Grafakos/Tao/Terwilleger. It was also extended to larger function spaces, first by Hunt to $L_p$ for $p > 1$ or even $L (\log L)^2$ and then by Antonov to $L \log L \log \log \log L$. Moreover, generalizations in terms of polynomial phase or the $p$-variation of partial sums have recently appeared in the literature, see \cite{Ant, Fef, GTT, Hun, LT, Lie, Taoetal} and the references therein. It is however worth mentioning that, despite so much work in this direction, very little is known on convergence of Fourier series for vector-valued functions. In 1986, J.L. Rubio de Francia \cite{RdF} proved that Carleson convergence theorem still holds for vector-valued functions taking values in a UMD Banach lattice. His argument uses ultimately the lattice structure to prove the result \lq point by point\rq${}$ and thereby reduce it to Carleson's statement. In the light of it, Rubio de Francia conjectured that Carleson theorem should hold for all UMD Banach spaces. Particularly, he pointed the Schatten $p$-classes for $1 < p < \infty$ as the simplest models for UMD Banach spaces not being a lattice. Apparently, no progress has been made since then. In this note we provide a step forward by proving the \lq little Carleson theorem\rq${}$ on the growth of Fourier series for arbitrary UMD Banach spaces. 

\begin{TheoA} If $f \in L (\log L)^{1 + \delta}(\T; \mathrm{X})$, we have $$\Big\| \sum_{|k| \le n} \widehat{f}(k) e^{2\pi i k x} \Big\|_\mathrm{X} \, = \, o \big( \log \log n \big) \quad \mbox{for} \quad \mathrm{ae} \hskip1pt - \hskip1pt x \in \T$$ provided $\delta > 0$ and $\mathrm{X}$ is a Banach space with unconditional martingale differences.
\end{TheoA}

In contrast to Rubio de Francia's approach, our proof is modeled by the full strength of Carleson's original argument to include nonlattice UMD spaces. Beyond some standard modifications needed in the vector-valued setting, it just presents one substantial departure from the scalar-valued argument. Specifically,  Zygmund's map in \cite[page 158]{Z} does not lead anymore to the corresponding Hausdorff-Young type inequality in $L (\log L)^{1+\delta}$ for UMD spaces. This estimate is crucial to control the size of the exceptional set in Carleson's approach. Here we modify Zygmund's construction to make it work on Banach spaces with nontrivial Fourier type, a condition which is even less restrictive that the UMD property. The resulting inequality is of independent interest. In what follows, se shall use $- \hskip-9pt \int_w$ to denote the mean $\frac{1}{|w|} \int_w$ over a finite interval $w$.

\begin{TheoB}
Assume that $$\mean_w \|f(x)\|_\mathrm{X} \big( \log^+ \|f(x)\|_\mathrm{X} \big)^{\beta} \, dx \, \le \, \rho$$ over a finite interval $w$ in $\R$ for some $\beta > 0$. Then, we have $$\sum_{k \in \Z} \exp \Big( - a(\rho,\beta) \|\widehat{f}_w(k)\|_\mathrm{X}^{- \frac{1}{\beta}} \Big) \ \le \ A(\rho,\beta)$$ for some constants $a(\rho,\beta), A(\rho,\beta) > 0$, provided $\mathrm{X}$ has non trivial Fourier type.
\end{TheoB}

We have decided to present a self-contained proof of Theorem A ---despite the parallelism with Carleson's original argument--- which we believe will help the reader who is not familiar with Carleson's paper. As it is well-known, we may replace the truncated Fourier series in the statement of Theorem A by the so-called \emph{modified partial sums} $$S_nf(x) = \int_\T \frac{f(t) e^{- 2 \pi i nt}}{x-t} \, dt.$$ We claim that for any $\varepsilon > 0$, there exists a measurable $\Sigma_{f, \varepsilon} \subset \T$ such that
\begin{itemize}
\item $|\Sigma_{f,\varepsilon}| < \varepsilon$,

\vskip3pt

\item $\big\| S_nf(x) \big\|_\mathrm{X} \le M(f,\varepsilon) \log \log n$ for all $x \notin \Sigma_{f, \varepsilon}$,

\vskip3pt

\item $M(f,\varepsilon)$ is arbitrarily small for $\|f\|_{L (\log L)^{1 + \delta}}$ small and $\varepsilon$ fixed.
\end{itemize}
It is apparent that our result follows from the claim above, since the $\mathrm{X}$-valued trigonometric polynomials form a dense subspace of $L (\log L)^{1+\delta}(\T;\mathrm{X})$ for which our result holds trivially. Section \ref{Sect1} contains some well-known preliminary estimates for UMD Banach spaces. In Section \ref{Sect2} we give the proof of Theorem B. In the rest of the paper we present Carleson decomposition of $f$, construct the exceptional sets $\Sigma_{f,\varepsilon}$ and finally complete the proof of Theorem A.

\section{Maximal Hilbert transform}
\label{Sect1}

The characterization of UMD Banach spaces in terms of the $L_p$ boundedness of vector-valued Calder\'on-Zygmund operators goes back to results by Bourgain and Burkholder in the 1980s, further generalized by Figiel. We refer to \cite{Bur2} for a nice survey paper on this subject. Let $$Hf(x) = \mathrm{p.v.} \int_\R \frac{f(y)}{x-y} \, dy,$$ the Hilbert transform acting on a vector-valued function $f: \R \to \mathrm{X}$. Given a finite interval $w$ in $\R$ and $x \in w$, let $\mathcal{I}_w(x)$ stand for the family of all subintervals $\gamma \subset w$ such that $x \in \frac12 \gamma$. We will also deal with the maximal Hilbert transform $$H_w^*f(x) = \sup_{\gamma \in \mathcal{I}_w(x)} \Big\| \mathrm{p.v} \int_\gamma \frac{f(y)}{x-y} \, dy \Big\|_\mathrm{X}.$$ 

\begin{remark} \label{RemH1H2}
\emph{The usual definition of the maximal Hilbert transform is $$H^*f(x) = \sup_{\varepsilon > 0} \Big\| \int_{|x-y| > \varepsilon} \frac{f(y)}{x-y} \, dy \Big\|_\mathrm{X}.$$ The given definition was introduced by Carleson and it is adapted to the proof of Theorem A. Note however that both operators are comparable. Indeed, given $x \in \R$ and an interval $\gamma$ such that $x \in \frac12 \gamma$, let $\sigma_x \subset \gamma$ denote the interval centered at $x$ of maximal length. Then, we find that $$\Big\| \mathrm{p.v} \int_\gamma \frac{f(y)}{x-y} \, dy \Big\|_\mathrm{X} \, \le \, \big\| Hf(x) \big\|_\mathrm{X} + \Big\| \int_{\R \setminus \sigma_x} \frac{f(y)}{x-y} \, dy \Big\|_\mathrm{X} + \Big\| \int_{\gamma \setminus \sigma_x} \frac{f(y)}{x-y} \, dy \Big\|_\mathrm{X}.$$ Moreover, if we set $$Mf(x) = \sup_{x \in \gamma} \mean_\gamma \|f(y)\|_\mathrm{X} \, dy,$$ and recall that $|x-y| \ge \frac14 |\gamma|$ whenever $y \in \gamma \setminus \sigma_x$ and $x \in \frac12 \gamma$, we conclude $$H_w^*f(x) \, \le \, \|Hf(x)\|_\mathrm{X} + H^*f(x) + 4 Mf(x).$$}
\end{remark}

\noindent In the following lemma we outline some well-known estimates to be used below.

\begin{lemma} \label{MaximalH}
If $\mathrm{X}$ is a \emph{UMD} Banach space, then $H: L_p(\R;\mathrm{X}) \to L_p(\R;\mathrm{X})$ is bounded for $1 < p < \infty$. Moreover, its norm is dominated by $c_\mathrm{X} \, p^2/(p-1)$ where $c_\mathrm{X}$ only depends on the \emph{UMD} constant of $\mathrm{X}$. In particular, we find the following estimate for some $c_0 > 0$ $$|w| \, \Big| \Big\{ x \in w \, \big| \ H^*_wf(x) > \lambda \Big\} \Big| \, \lesssim \, \min \Big\{ \frac{1}{\lambda} \, \mean_w \|f(x)\|_\mathrm{X} \, dx, \ \exp \big( - c_0 \lambda/\|f\|_\infty \big) \Big\}.$$
\end{lemma}

\ske According to \cite{Bou,Bur}, the UMD property is equivalent to the boundedness of $H: L_p(\R;\mathrm{X}) \to L_p(\R;\mathrm{X})$. Once we know the boundedness on the space (say) $L_2(\R; \mathrm{X})$, the behavior of the constants follows by interpolation with the usual endpoint spaces, since the classical arguments apply here to prove the weak type $(1,1)$ and the $L_\infty \to \mathrm{BMO}$ estimates. Both estimates for the maximal Hilbert transform rely on Cotlar's inequality $$H^*f(x) \, \lesssim \, MHf(x) + Mf(x),$$ which holds for arbitrary Banach spaces. According to Remark \ref{RemH1H2}, since the weak type $(1,1)$ inequality holds for the Hilbert transform and the maximal function, our first estimate reduces to the same one for $H^*f$. A standard regularization argument yields $$\big| H^*f(x) - H_\phi^*f(x) \big| \, \lesssim \, Mf(x),$$ where $H_\phi^*f(x)$ is defined replacing $\chi_{|x-y|> \varepsilon}$ in $H^*f(x)$ by $\phi(\frac{1}{\varepsilon} |x-y|)$ with $\phi$ a smooth function which vanishes on $(-\infty, 0.5]$ and is identically $1$ on $[1.5, \infty)$. Now it can be checked that $$H_\phi^*f(x) \, = \, \Big\| \int_\R K_\phi(x-y) f(y) \, dy \Big\|_{\ell_\infty(\R_+; \mathrm{X})} \quad \mbox{with} \quad K_\phi(x) = \Big\{ \frac{1}{x} \phi \big(\frac{|x|}{\varepsilon} \big) \Big\}_{\varepsilon >0}$$ a Calder\'on-Zygmund kernel satisfying the usual Lipschitz estimates. Hence, the weak type estimate will follow from a Calder\'on-Zygmund decomposition if we know that $T_\phi: L_2(\R; \mathrm{X}) \to L_2(\R; \ell_\infty(\mathrm{X}))$ is bounded, which in turn follows from Cotlar's inequality above. Let us now prove the exponential type estimate. By homogeneity we may assume that $\|f\|_\infty = 1$. On the other hand, Cotlar's inequality gives in conjunction with Remark \ref{RemH1H2} that $$H_w^*f(x) \, \lesssim \, \|Hf(x)\|_\mathrm{X} + Mf(x) + MHf(x).$$ Therefore, according to the first part of the statement we see that 
\begin{eqnarray*}
\lefteqn{\hskip-10pt \Big| \Big\{ x \in w \, \big| \ H^*_wf(x) > \lambda \Big\} \Big|} \\ & \lesssim & \min \Big\{ |w|, \frac{c_\mathrm{X}^n n^n}{\lambda^n} \|f \chi_w\|_n^n \, \big| \, n \ge 3  \Big\} \ \le \ |w| \min \Big\{ 1, \frac{c_\mathrm{X}^n n^n}{\lambda^n} \, \big| \, n \ge 3 \Big\}.
\end{eqnarray*}
Since $e^x \, \lesssim \, 1 + \sum_{n \ge 3} x^n/n!$ up to an absolute constant, we find $$\Big| \Big\{ x \in w \, \big| \ H^*_wf(x) > \lambda \Big\} \Big| \, \lesssim \, e^{- c_0 \lambda} |w| \Big( 1 + \sum_{n \ge 3} \frac{c_0^n c_\mathrm{X}^n n^n}{n!} \Big) \, \lesssim \, e^{- c_0 \lambda} |w|.$$ The last inequality follows from Stirling's formula by picking $c_0$ small enough. \fin

\section{A Hausdorff-Young type inequality}
\label{Sect2}

The UMD condition is a super property. Since $\ell_1$ is not a UMD space, all UMD spaces fail to contain $\ell_1^n$'s uniformly, which is known to be equivalent to having nontrivial type. According to Bourgain \cite{Bou2}, we find that every UMD Banach space satisfies a nontrivial Hausdorff-Young inequality. In other words, given a UMD Banach space $\mathrm{X}$ there exists some $1 < p \le 2$ such that $$\Big( \sum_{k \in \Z} \|\widehat{f}(k)\|_\mathrm{X}^q \Big)^{\frac1q} \, \le \, c_p \Big( \int_\T \|f(x)\|_\mathrm{X}^p \, dx \Big)^{\frac1p} \quad \mbox{with} \quad \frac1p + \frac1q = 1.$$ A Banach space satisfying this inequality is said to have Fourier type $p$. Note that Fourier type $1$ trivially holds for every Banach space. In this section we prove Theorem B, a Hausdorff-Young type inequality on $L (\log L)^\beta$ for Banach spaces with non trivial Fourier type. Given a finite interval $w$ in $\R$, let us equip it with its normalized measure $d\mu(t) = dt / |w|$. Let us also fix $\alpha > 1$ and consider the measure on $\Z$ given by $\nu_\alpha (\{k\}) = |k|^{-\alpha}$, except for $k=0$ where we impose $\nu_\alpha (\{0\}) = 1$. Given a permutation $\sigma: \Z \to \Z$, define the linear map $$\Lambda_{\sigma}f = \big( |k|^{\alpha-1} \widehat{f}_w(\sigma(k)) \big)_{k \in \Z} \quad \mbox{with} \quad \widehat{f}_w(k) = \mean_w f(t) e^{- 2 \pi i k t / |w|} \, dt.$$ 

\begin{lemma} \label{LemmaHY1}
If $1 < \alpha < 2$, $\Lambda_{\sigma}$ satisfies
\begin{itemize}
\item[i)] $\Lambda_{\sigma}: L_1(w, \mu; \mathrm{X}) \to L_{1,\infty}(\Z, \nu_\alpha; \mathrm{X})$ is bounded for any $\mathrm{X}$.

\item[ii)] If $\mathrm{X}$ has Fourier type $p >1$, $\Lambda_{\sigma}: L_p(w, \mu; \mathrm{X}) \to L_p(\Z,\nu_\alpha; \mathrm{X})$.
\end{itemize}
\end{lemma}

\dem The first assertion follows from $$\sum_{\begin{subarray}{c}k \in \Z \\ |k|^{\alpha-1} \|\widehat{f}_w(\sigma(k))\|_\mathrm{X} > \lambda \end{subarray}} |k|^{-\alpha} \, \le \, \frac{c_\alpha}{\lambda} \|f\|_1.$$ Indeed, the summation index is contained in the set of integers $$|k| > \Big( \frac{\lambda}{\|f\|_1}\Big)^{\frac{1}{\alpha-1}}, \quad \mbox{since} \quad \|\widehat{f}_w(\sigma(k))\|_\mathrm{X} \le \mean_w \|f(t)\|_\mathrm{X} \, dt = \|f\|_1.$$ If $\mathrm{X}$ has Fourier type $p > 1$, we have for $\frac1p + \frac1q = 1$
\begin{eqnarray*}
\|\Lambda_{\sigma} f\|_p & = & \Big( \sum_{k \in \Z \setminus \{0\}} \frac{\big\| |k|^{\alpha-1} \widehat{f}_w(\sigma(k)) \big\|_\mathrm{X}^p}{|k|^{\alpha}} \Big)^{\frac{1}{p}} \\ & \le & \Big( \sum_{k \in \Z \setminus \{0\}} \Big[ \frac{1}{|k|^{\alpha - (\alpha-1)p}} \Big]^{\frac{q}{q-p}} \Big)^{\frac{q-p}{pq}} \Big( \sum_{k \in \Z} \|\widehat{f}_w(k) \|_\mathrm{X}^q \Big)^{\frac{1}{q}} \\ & \le & c_{\mathrm{X},\alpha} \Big( \mean_w \|f(t)\|_\mathrm{X}^p \, dt \Big)^{\frac{1}{p}} \ = \ c_{\mathrm{X},\alpha} \, \|f\|_p.
\end{eqnarray*}
The last inequality uses $\frac{q[\alpha - (\alpha-1)p]}{q-p} > 1$ iff $\alpha < 2$. The proof is complete. \fin

\begin{remark}
\emph{Lemma \ref{LemmaHY1} yields in fact a characterization of Banach spaces with nontrivial Fourier type. In other words, $\mathrm{X}$ has nontrivial Fourier type iff $\Lambda_\sigma$ is $L_p(\mathrm{X})$-bounded for some $p > 1$. The sufficiency appears in the proof above. For the necessity, it suffices to show that $\mathrm{X} = L_1(\T)$ fails inequality ii) in the Lemma for any $p > 1$. Indeed, since our inequality extends to finitely representable spaces in $\mathrm{X}$, our claim for $L_1(\T)$ means that if $\Lambda_\sigma$ is $L_p(\mathrm{X})$-bounded, then $\mathrm{X}$ can not contain $\ell_1^n$'s uniformly, which in turn characterizes nontrivial Fourier type. The counterexample arises from the Poisson kernel, take $f_r: \T \to L_1(\T)$ given by $$f_r(x) = P_r(x + \cdot \hskip1pt ) = \sum_{k \in \Z} ( r^{|k|} e^{2\pi i k \cdot} ) e^{2\pi i k x} = \sum_{k \in \Z} \widehat{f}(k) e^{2\pi i k x}.$$ Since $\|f_r(x)\|_{L_1(\T)} = 1$ for all $x \in \T$ and $0 < r < 1$, we get easily the conclusion.}
\end{remark}

\begin{lemma} \label{LemmaHY2}
Given $\beta > 0$, we have for $f \in L (\log L)^\beta(w; \mathrm{X})$ $$\sum_{|k| \ge 2} \frac{\|\widehat{f}_w(\sigma(k))\|_\mathrm{X}}{|k|} (\log |k|)^{\beta-1} \, \lesssim \, 1 + \mean_w \|f(t)\|_\mathrm{X} (\log^+ \|f(t)\|_\mathrm{X})^{\beta} \, dt.$$
\end{lemma}

\dem We claim that $$\sum_{k \in \Z \setminus \{0\}} \frac{|k|^{\alpha-1} \|\widehat{f}_w(\sigma(k))\|_\mathrm{X}}{|k|^\alpha} \big[ \log^+ \big( |k|^{\alpha-1} \|\widehat{f}_w(\sigma(k))\|_\mathrm{X} \big) \big]^{\beta-1}$$ is bounded above by the right hand side of the stated inequality. Indeed, this follows from Lemma \ref{LemmaHY1} together with \cite[Theorem 4.34, pag 118]{Z} applied to $\Lambda_{\sigma}$ with $\chi(u) = u (\log^+ u)^{\beta-1}$ and $\phi(u) = u (\log^+ u)^{\beta}$. Now, fixing the value of $\alpha = \frac{3}{2}$, the terms in the left hand side of the statement satisfying $$|k|^{-\frac{1}{4}} \le \|\widehat{f}_w(\sigma(k))\|_\mathrm{X} \le |k|^{\frac{1}{4}}$$ are comparable to the corresponding terms in the sum above. This completes our estimate for the main part of the sum. The terms satisfying $\|\widehat{f}_w(\sigma)(k)\|_\mathrm{X} < |k|^{- \frac{1}{4}}$ are bounded by $\sum_{|k| \ge 2} (\log |k|)^{\beta - 1} |k|^{- (1 + \frac{1}{4})} \lesssim 1$. Finally, since $$\|\widehat{f}_w(\sigma(k))\|_\mathrm{X} \le \|f\|_1,$$ the remaining terms satisfy $2 \le |k| \le \|f\|_1^4$ and the sum is dominated by $$\|f\|_1 \sum_{2 \le |k| \le \|f\|_1^4} \frac{(\log |k|)^{\beta-1}}{|k|} \le \|f\|_1 \big( \log^+ \|f\|_1 \big)^{\beta} \le \|f\|_{L (\log L)^{1+\delta}}.$$ The last estimate follows from Jensen inequality. This completes the proof. \fin

\demB Pick $\sigma$ so that $$\widehat{f}_w(\sigma(0)) \ge \widehat{f}_w(\sigma(1)) \ge \widehat{f}_w(\sigma(-1)) \ge \widehat{f}_w(\sigma(2)) \ge \widehat{f}_w(\sigma(-2))...$$ This and Lemma \ref{LemmaHY2} yield for $k$ positive $$\|\widehat{f}_w(\sigma(k))\|_\mathrm{X} \, (\log k)^\beta \, \lesssim \, \sum_{j= \sqrt{k}}^k \frac{\|\widehat{f}_w(\sigma(j))\|_\mathrm{X}}{k} (\log j)^{\beta - 1} \rightarrow 0 \quad \mbox{as $k \to \infty$}.$$ A similar argument for $k$ negative leads us to the conclusion that $$\widehat{f}_w^*(k) = o \big( (\log |k|)^{- \beta} \big),$$ where $\widehat{f}_w^*(k)$ stands for the decreasing rearrangement of the sequence of Fourier coefficients. Since the constants so far only depend on $\rho$ and $\beta$, there exists a universal $k_0 = k_0(\rho,\beta)$ such that $\|\widehat{f}_w^*(k)\|_\mathrm{X} (\log |k|)^{\beta} \le 1$ for all $|k| \ge k_0$. Let us consider the constant $$a(\rho,\beta) = \max \Big\{ 2, \rho^{\frac{1}{\beta}} \log k_0(\rho,\beta) \Big\}.$$ Then, we find 
\begin{eqnarray*}
\sum_{k \in \Z} \exp \Big( - a \|\widehat{f}_w(k)\|_\mathrm{X}^{-\frac{1}{\beta}} \Big) & = & \sum_{k \in \Z} \exp \Big( - a \|\widehat{f}_w^*(k)\|_\mathrm{X}^{-\frac{1}{\beta}} \Big) \\ & \le & \sum_{|k| < k_0} \exp \Big( - \rho^{\frac{1}{\beta}} \log k_0 \, \|\widehat{f}_w^*(k)\|_\mathrm{X}^{-\frac{1}{\beta}} \Big) \\ & + & \sum_{|k| \ge k_0} \exp \Big( -2 \big[ \|\widehat{f}_w^*(k)\|_\mathrm{X} (\log |k|)^{\beta} \big]^{-\frac{1}{\beta}} \log |k| \Big). 
\end{eqnarray*}
Since $\|\widehat{f}_w^*(k)\|_\mathrm{X} \le - \hskip-9pt \int_w \|f\|_\mathrm{X} \le \rho$, this sum is dominated by $A = 2 + \sum_{|k| \ge k_0} \frac{1}{k^2}$. \fin

\section{Carleson decomposition}
\label{Sect3}

Now we are ready to start the proof of Theorem A. In this section, we describe Carleson decomposition of $f$. This requires to introduce some terminology. In what follows, we will represent $\T$ by the interval $(-\frac12, \frac12)$, $w$ and $w'$ will denote dyadic intervals in $\T$ and $w_{-1}^*$ will stand for the interval $(-2,2)$. Let us set 

\begin{itemize}
\item Dyadic intervals in $(-1,1)$ $$\Big\{ w_{j\nu} \subset (-1,1) \ \big| \ |w_{j\nu}| = 2^{-\nu}, \ \nu \ge 0, \ 1 \le j \le 2^{\nu + 1} \Big\}.$$

\item Smoothing intervals, $w_{j\nu}^* = w_{j \nu} \cup w_{j+1, \nu}$ with $1 \le j \le 2^{\nu+1}-1$.

\vskip5pt

\item Generalized Fourier coefficients $$\widehat{f}_w(\alpha) = \mean_w f(t) e^{-2 \pi i \alpha t/|w|} \,
dt \quad \mbox{for} \quad \alpha \in \R.$$ Note that $\{e^{2 \pi i k \cdot / |w|}\}_{k \in \Z}$ forms
an orthonormal basis of $L_2(w, dt/|w|)$.

\vskip5pt

\item Carleson averages of Fourier coefficients $$\qquad C_k(f,w) = \frac{1}{\gamma}
\sum_{\mu \in \Z} \frac{\|\widehat{f}_w(k+\mu/3)\|_\mathrm{X}}{1 + \mu^2} \quad \mbox{with} \quad \gamma = \sum_{\mu \in \Z} \frac{1}{1+\mu^2}$$ and $k \in \Z$. \hskip-1pt We have $C_k(f,w) \le - \hskip-9pt \int_w
\|f\|_\mathrm{X}$ and $C_k(f,w)=0$ iff $f=0$ ae-$w$.

\vskip5pt

\item Amplified averages $$\qquad C_k^*(f,w^*) = \max \Big\{ C_k(f,w') \, \big| \, w' \subset w^* \ \mbox{and} \ 4|w'|=w^* \Big\}.$$ That is, we consider the Carleson averages for the dyadic grandsons of $w^*$. 
\end{itemize}

Given $k \in \N$, we also set $k[w] = [k|w|]$ ---the integer part of $k|w|$--- for dyadic intervals and $k[w^*] = [\frac14 k|w^*|]$ for smoothing intervals. Given a smoothing interval $w^*$, a nonnegative integer $k$ and $\lambda > 0$, assume $$C_{k[w^*]}^*(f,w^*) \, \le \, \lambda.$$ We will write $\Omega_\lambda(k,w^*)$ for the corresponding Carleson partition of $w^*$, which is constructed as follows. Each element of our partition will be a proper dyadic subinterval $w'$ of $w^*$ satisfying $|w'| \ge 1/2n$ and
\begin{itemize}
\item[i)] $C_{k[w']}(f,w') \le \lambda$.
\end{itemize}
These conditions however do not determine a unique partition. For instance, the 4 dyadic grandsons of $w^*$ satisfy i) when $|w^*| \ge 2/n$. The additional conditions to impose our choice to be maximal are as follows
\begin{itemize}
\item[ii)]A dyadic son of $w'$ fails i) or $|w'| = 2^{-[\log_2 2n]-1}$,

\vskip3pt

\item[iii)] $w'$ is maximal among the intervals satisfying i) and ii). 
\end{itemize}
We will use $\Omega_\lambda(k,w^*)$ to decompose the function $f_k(x) = f(x) \exp(- 2 \pi i k x) \chi_{w^*}(x)$ as $f_k(x) = g_{k,\lambda}(x) + b_{k,\lambda}(x)$, where the \lq good/bad\rq${}$ parts $g_{k,\lambda}$ and $b_{k,\lambda}$ are given by $$g_{k,\lambda} = \sum_{w' \in \Omega_\lambda(k,w^*)} \Big( \mean_{w'} f(t) \exp(-2\pi i k t) \, dt \Big) \, \chi_{w'}(x) = \sum_{w' \in \Omega_\lambda(k,w^*)} f_{k}[w'] \chi_{w'}(x)$$ and $b_{k,\lambda}(x) = f_k(x) - g_{k,\lambda}(x) = \sum_{\Omega_\lambda(k,w^*)} (f_k(x) - f_k[w']) \chi_{w'}(x)$. We will refer to this as Carleson decomposition. A moment of thought shows many similarities between Carleson and Calder\'on-Zygmund decompositions. Indeed, let us consider the following maximal function $$\mathcal{M}_kf(x) = \sup \Big\{ C_{k[w]}(f,w) \, \big| \, x \in w \subset w^* \Big\},$$ the analog of the dyadic maximal function defined from Carleson averages instead of dyadic ones. Ignoring the size truncation $|w'| \ge 1/2n$, $\Omega_\lambda(k,w^*)$ would be the union of dyadic fathers of maximal intervals for $\{\mathcal{M}_kf > \lambda\}$. Then, up to this shifted generation, Carleson decomposition also follows the usual averaging/deaveraging procedure which we find in Calder\'on-Zygmund decomposition. Carleson averages (instead of usual ones) are crucial to estimate the size of the exceptional set.

\begin{lemma} \label{Lemmagoodpart} 
In the situation above, we have $\|g_{k,\lambda}\|_\infty \lesssim \lambda$.
\end{lemma}

\dem It suffices to prove the estimate $$\|f_k[w']\|_\mathrm{X} = \|\widehat{f}_{w'}(k|w'|)\|_\mathrm{X} \lesssim C_{k[w']}(f,w'),$$ since the right hand side is bounded by $\lambda$ due to the construction of $\Omega_\lambda(k,w^*)$. If $k|w'|$ is an integer, the inequality is clear since the left hand side is the $\mu=0$ term on the right. Otherwise, 
we write $\exp(-2\pi i k x) = \exp(-2 \pi i \beta x) \exp(-2\pi i k[w'] x / |w'|)$ for $\beta = k - \frac{1}{|w'|}k[w']$ and expand $$\exp(- 2 \pi i \beta x) = \sum_{\mu \in \Z} \gamma_\mu \exp(-2 \pi i \mu x / 3|w'|) \quad \mbox{for} \quad x \in w'.$$ Namely, extend $\phi_\beta(x) = \exp(-2 \pi i \beta x) \chi_{w'}(x)$ to a smooth, compactly supported function in $3w'$. The expression above then follows as the Fourier series adapted to $3w'$. Integration by parts gives $$(1+\mu^2) |\gamma_\mu| \, \lesssim \, \|\phi_\beta\|_\infty + |w'|^2 \|\phi_\beta''\|_\infty \, \lesssim \, 1$$ and the result follows. This argument appears in \cite[Lemmas 2 and 3]{C}. \fin

\vskip5pt

\section{The exceptional set}
\label{Sect4}

In this section we construct the exceptional set $\Sigma_{f,\varepsilon}$ and estimate its size. $\Sigma_{f,\varepsilon}$ will be the union of four sets $\Sigma_{f,\varepsilon}^j$, $1 \le j \le 4$. Since we are assuming in Theorem A that $f \in L (\log L)^{1 + \delta}(\T; \mathrm{X})$, we set $$\lambda = \frac{\|f\|_{L(\log L)^{1+ \delta}}}{\varepsilon}.$$ Then, if $\varphi(x) = \|f(x)\|_\mathrm{X} (\log^+ \|f(x)\|_\mathrm{X})^{1+\delta}$, we define $$\Sigma_{f,\varepsilon}^1 = 7 \Big\{ M_d \varphi > \lambda \Big\} \quad \mbox{and} \quad \Sigma_{f,\varepsilon}^2 = \Big\{ H_\T^*f > \lambda \Big\}.$$ Here $M_d$ stands for the dyadic Hardy-Littlewood maximal operator and the factor 7 means that we dilate concentrically each maximal interval in $\{M_d \varphi > \lambda\}$ by this factor. To define the other pieces of the exceptional set, we consider a pair $(k,w^*)$ such that $\frac14 k|w^*| \in \Z$. If $w^* \not\subset \Sigma_{f,\varepsilon}^1$, no grandson of $w^*$ may belong to the $\lambda$-level set of $M_d \varphi$ and we get $$C_{k[w^*]}^*(f,w^*) \, = \, \max_{\mathrm{grandsons}} C_{k[w']}(f,w') \, \le \, \max_{\mathrm{grandsons}} \mean_{w'} \|f(x)\|_\mathrm{X} \, dx \, \le \, \lambda.$$ Thus, for any such $(k,w^*)$ there exists a unique $j \ge 1$ such that $$2^{-j} \lambda \, < \, C_{k[w^*]}^*(f,w^*) \, \le \, 2^{1-j} \lambda,$$ unless $f \equiv 0$ over $w^*$. These $j$'s allow us to introduce the sets 
\begin{eqnarray*}
\Sigma_{f,\varepsilon}^3(k,w^*,n) & = & \Big\{  H_{w^*}^*g_{k,2^{1-j} \lambda} \hskip1pt > \hskip1pt R(\lambda) (2^{1-j} \lambda)^{1- \frac{1}{1+\delta/2}} \log\log n \Big\}, \\ \Sigma_{f,\varepsilon}^4(k,w^*,n) & = & \Big\{ \Delta_{\Omega_{2^{1-j} \lambda}(k,w^*)} > R(\lambda) (2^{1-j} \lambda)^{- \frac{1}{1+\delta/2}} \log\log n \Big\}.
\end{eqnarray*}
Here, $R(\lambda)$ is a constant to be fixed and $$\Delta_{\Omega_{2^{1-j}\lambda}(k,w^*)}(x) = \sum_{w' \in \Omega_{2^{1-j}\lambda}(k,w^*)} \frac{|w'|^2}{(x - \mathrm{c}(w'))^2 + |w'|^2}$$ with $\mathrm{c}(w')$ the center of $w'$. Consider the family of pairs $$\mathcal{A}_n \, = \, \Big\{ (k,w^*) \, \big| \, w^* \not\subset \Sigma_{f,\varepsilon}^1, \ \frac14 k|w^*| \in \Z, \ 1 \le k \le n \ \mbox{and} \ |w^*| \ge 2^{-[\log_2 2n]} \Big\}.$$ Then, we may define the other pieces of the exceptional set as follows $$\Sigma_{f,\varepsilon}^j \, = \, \bigcup_{m \ge \frac{1}{\varepsilon} C(\lambda,\delta)} \bigcup_{(k,w^*) \in \mathcal{A}_{e^m}} \Sigma_{f,\varepsilon}^j(k,w^*,e^m) \qquad \mbox{for} \qquad j=3,4,$$ where the constant $C(\lambda,\delta)$ is also to be fixed. Our goal in the rest of this section is to estimate the size of the exceptional set $\Sigma_{f,\varepsilon} = \Sigma_{f,\varepsilon}^1 \cup \Sigma_{f,\varepsilon}^2 \cup \Sigma_{f,\varepsilon}^3 \cup \Sigma_{f,\varepsilon}^4$. A key result is the following. 

\begin{lemma} \label{LemmaHY3} 
Assume that $$\mean_w \|f(x)\|_\mathrm{X} \big( \log^+ \|f(x)\|_\mathrm{X} \big)^{1+\delta} \, dx \, \le \, \lambda.$$ Then, there exist $b(\lambda,\delta), B(\lambda,\delta) > 0$ such that $$\sum_{k \in \Z} \exp \Big( - b(\lambda, \delta) C_k(f,w)^{- \frac{1}{1+\delta}} \Big) \ \le \ B(\lambda,\delta).$$
\end{lemma}

\dem We just need to follow the argument \cite[Lemma 1]{C} replacing Zygmund's results by our results from Section \ref{Sect2}. According to Theorem B, the statement holds for $\widehat{f}_w(k)$ instead of $C_k(f,w)$. After modulating $f$ with $\exp(\pm \frac23 \pi i x)$, we see that Theorem B also holds for frequencies $$\widehat{f}_w(k \pm \frac13).$$ Since $$C_k(f,w) \, \lesssim \, \sup_{\mu \in \Z} \frac{\|\widehat{f}_w(k+\mu/3)\|_\mathrm{X}}{\sqrt{1 + |\mu|}}$$ and $\|\widehat{f}_w(k+\mu/3)\|_\mathrm{X} \le - \hskip-9pt \int_w \|f\|_\mathrm{X} \le \lambda$, we find 
\begin{eqnarray*}
\lefteqn{\sum_{k \in \Z} \exp \Big( - b(\lambda, \delta) C_k(f,w)^{- \frac{1}{1+\delta}} \Big)} \\ & \le & \sum_{k \in \Z} \sum_{\mu \in \Z} \exp \Big( - a (\lambda, \delta) \|\widehat{f}_w(k+\mu/3)\|_\mathrm{X}^{- \frac{1}{1+\delta}} \sqrt{1 + |\mu|}^{\frac{1}{1+\delta}} \Big) \\ & \lesssim & \sum_{\mu \in \Z} \sum_{k \in \Z} \exp \Big( - a (\lambda, \delta) \|\widehat{f}_w(k+\mu/3)\|_\mathrm{X}^{- \frac{1}{1+\delta}} \Big) \exp \Big( - a (\lambda, \delta) \lambda^{- \frac{1}{1+\delta}} \sqrt{|\mu|}^{\frac{1}{1+\delta}} \Big) \\ & \le & 3 A(\lambda, \delta) \sum_{\mu \in \Z} \exp \Big( - a(\lambda, \delta) \lambda^{- \frac{1}{1+\delta}} \sqrt{|\mu|}^{\frac{1}{1+\delta}} \Big) \ \le \ B(\lambda,\delta)
\end{eqnarray*}
for some constant $b(\lambda,\delta) \sim a(\lambda, \delta)$. This completes the proof. \fin 

\begin{lemma} Continuing with the exceptional set, we have
$$|\Sigma_{f,\varepsilon}| \, = \, \Big| \bigcup_{j=1}^4 \Sigma_{f,\varepsilon}^j \Big| \, \lesssim \, \varepsilon.$$
\end{lemma}

\dem Clearly $|\Sigma_{f,\varepsilon}^1| \le 7 \varepsilon$ and $$|\Sigma_{f,\varepsilon}^2| \, \le \, \frac{1}{\lambda} \int_\T \|f(x)\|_\mathrm{X} \, dx \, \le \, \varepsilon$$ from Lemma \ref{MaximalH}. To estimate $|\Sigma_{f,\varepsilon}^3|$ and $|\Sigma_{f,\varepsilon}^4|$, we claim $$\sum_{(k,w^*) \in \mathcal{A}_n} |\Sigma_{f,\varepsilon}^3(k,w^*,n)| + |\Sigma_{f,\varepsilon}^4(k,w^*,n)| \, \lesssim \, C(\lambda,\delta) \frac{1}{(\log n)^2}.$$ It is clear that the statement follows from the claim, since we have for $j=3,4$ $$|\Sigma_{f,\varepsilon}^j| \, \le \, \sum_{m \ge \frac{1}{\varepsilon} C(\lambda,\delta)} \sum_{(k,w^*) \in \mathcal{A}_{e^m}} |\Sigma_{f,\varepsilon}^j(k,w^*,e^m)| \, \le \, C(\lambda,\delta) \sum_{m \ge \frac{1}{\varepsilon} C(\lambda,\delta)} \frac{1}{m^2} \, \lesssim \, \varepsilon.$$ Let us begin by considering the claim for $\Sigma_{f,\varepsilon}^3$. According to Lemmas \ref{MaximalH} and \ref{Lemmagoodpart} 
\begin{eqnarray*}
|\Sigma_{f,\varepsilon}^3(k,w^*,n)| & \lesssim & \exp \Big( - c_0 R(\lambda) (2^{1-j} \lambda)^{1- \frac{1}{1+\delta/2}} \log\log n / \|g_{k,2^{1-j}\lambda}\|_\infty \Big) \, |w^*| \\ & \lesssim & \exp \Big( - c_0 R(\lambda) (2^{1-j} \lambda)^{- \frac{1}{1+\delta/2}} \log\log n \Big) \, |w^*|
\end{eqnarray*}
Since $w^* \not\subset \Sigma_{f,\varepsilon}^1$ for $w^* \in \mathcal{A}_n$, we get $$\mean_{w'} \|f\|_\mathrm{X} (\log^+\|f\|_\mathrm{X})^{1+\delta} \le \lambda$$ for every grandson $w'$ of $w^*$. Lemma \ref{LemmaHY3} then gives 
\begin{eqnarray*}
\lefteqn{\hskip-10pt \sum_{(k,w^*) \in \mathcal{A}_n} \exp \Big( - b(\lambda,\delta) C_{k[w^*]}^*(f,w^*)^{- \frac{1}{1+\delta}} \Big) |w^*|} \\ & \le & \sum_{\begin{subarray}{c} w^* \subset \T, w^* \not\subset \Sigma_{f,\varepsilon}^1 \\ |w^*| \ge 2^{-[\log_2 2n]} \end{subarray}} |w^*| \Big[ \sum_{m \in \Z} \exp \Big( - b(\lambda,\delta) C_m^*(f,w^*)^{- \frac{1}{1+\delta}} \Big) \Big] \ \lesssim \ B(\lambda,\delta) \log n.
\end{eqnarray*}
We are now ready to prove the claim, define the sets $$\mathcal{A}_{nj} \, = \, \Big\{ (k,w^*) \in \mathcal{A}_n \, | \, 2^{-j} \lambda < C_{k[w^*]}^*(f,w^*) \le 2^{1-j} \lambda \Big\}.$$ Our estimates so far give rise to the following inequalities
\begin{eqnarray*}
\sum_{(k,w^*) \in \mathcal{A}_{nj}} |w^*| & \le & B(\lambda,\delta) \log n \, \exp \Big( b(\lambda,\delta) (2^{-j} \lambda)^{- \frac{1}{1+\delta}} \Big), \\ |\Sigma_{f,\varepsilon}^3(k,w^*,n)| & \lesssim & \exp \Big( - c_0 R(\lambda) (2^{1-j} \lambda)^{- \frac{1}{1+\delta/2}} \log\log n \Big) \, |w^*|.
\end{eqnarray*}
Then we may prove the claim for $\Sigma_{f,\varepsilon}^3$ using the partition $\mathcal{A}_n = \cup_{j \ge 1} \mathcal{A}_{nj}$ as follows
\begin{eqnarray*}
\lefteqn{\hskip-5pt \sum_{(k,w^*) \in \mathcal{A}_n} |\Sigma_{f,\varepsilon}^3(k,w^*,n)|} \\ & \lesssim & \sum_{j \ge 1} \Big[ \sum_{(k,w^*) \in \mathcal{A}_{nj}} |w^*| \Big] \exp \Big( - c_0 R(\lambda) (2^{1-j} \lambda)^{- \frac{1}{1+\delta/2}} \log\log n \Big) \\ & \lesssim & B(\lambda,\delta) \log n \sum_{j \ge 1} \exp \Big( b(\lambda,\delta) (2^{-j} \lambda)^{- \frac{1}{1+\delta}} - c_0 R(\lambda) (2^{1-j} \lambda)^{- \frac{1}{1+\delta/2}} \log\log n \Big)
\end{eqnarray*}
Fixing the value of $R(\lambda)$ as $$R(\lambda) \, = \, \frac{3}{c_0} \lambda^{\frac{1}{1 + \delta/2}} + \frac{b(\lambda,\delta)}{c_0} 2^{\frac{1}{1+\delta/2}} \lambda^{\frac{1}{1 + \delta/2} - \frac{1}{1+\delta}},$$ we obtain the following estimate 
\begin{eqnarray*}
\lefteqn{\hskip-20pt \sum_{(k,w^*) \in \mathcal{A}_n} |\Sigma_{f,\varepsilon}^3(k,w^*,n)|} \\ & \lesssim & \frac{B(\lambda,\delta)}{(\log n)^2} \sum_{j \ge 1} \exp \Big( b(\lambda,\delta) \lambda^{- \frac{1}{1+\delta}} (2^{\frac{j}{1+\delta}} - 2^{\frac{j}{1+\delta/2}}) \Big) \ \equiv \ C(\lambda,\delta) \frac{1}{(\log n)^2}.
\end{eqnarray*}
This completes the proof of the claim for $\Sigma_{f,\varepsilon}^3$. Note that the only point in the argument where we use specific information about $\Sigma_{f,\varepsilon}^3$ is our estimate for the measure $|\Sigma_{f,\varepsilon}^3(k,w^*,n)|$. Therefore, it suffices to show that the same estimate holds for $|\Sigma_{f,\varepsilon}^4(k,w^*,n)|$. In other words, we just need to prove that $$\Big| \Big\{ \Delta_{\Omega(w^*)} > \rho \Big\} \Big| \, \lesssim \, \exp \big( - c_0 \rho \big) |w^*|$$ for some absolute constant $c_0$, any $\rho > 0$ and every partition $\Omega(w^*)$ into disjoint intervals. Here we write $\Delta_{\Omega(w^*)}(x) = \sum_{w' \in \Omega(w^*)} |w'|^2/ ((x - \mathrm{c}(w'))^2 + |w'|^2)$ and the corresponding estimate is quite standard, see \cite[Lemma 5]{C}. \fin

\section{Growth of vector-valued Fourier series}
\label{Sect5}

In this section, we complete the proof of Theorem A. The goal is to show that $\|S_nf(x)\|_\mathrm{X} \le M(f,\varepsilon) \log \log n$ on $\T \setminus \Sigma_{f,\varepsilon}$, for some constant $M(f,\varepsilon) \to 0^+$ as $\|f\|_{L(\log L)^{1+\delta}} \to 0$ for $\varepsilon > 0$ fixed. Given $x \in \T \setminus \Sigma_{f,\varepsilon}$ and $n \ge 1$, we will construct a finite family of triplets $$\Big\{ (j_s, k_s,w_s^*) \, \big| \, -1 \le s \le L \Big\}$$ satisfying 
\begin{itemize}
\item[i)] $j_s \in \N$ and $j_s > j_{s+1}$,

\vskip3pt

\item[ii)] $x \in \frac12 w_s^*$ for all $s$, $w_{-1}^* = (-2,2)$ and $w_{s+1}^* \subsetneq w_s^*$,

\vskip3pt

\item[iii)] $\frac14 k_s |w_s^*| \in \Z$ and $0 = k_L \le \ldots \le k_{s+1} \le k_s \ldots \le k_{-1} = n$,

\vskip3pt

\item[iv)] There exists $M'(f,\varepsilon) \to 0^+$ as $f \to 0$ for fixed $\varepsilon$ such that
$$\hskip32pt \big\| S_{k_s}f(x,w_s^*) \big\|_{\mathrm{X}} = \big\| S_{k_{s+1}}f(x,w_{s+1}^*) \big\|_{\mathrm{X}} + O \big( 2^{-j_{s+1}(1 - \frac{1}{1+\frac{\delta}{2}})} M'(f,\varepsilon) \log \log n \big).$$
\end{itemize}
Here, $S_kf(x,w^*)$ are modified partial sums adapted to $w^*$ $$S_kf(x,w^*) \, = \, \int_{w^*} \frac{f(t) e^{-2 \pi i k t}}{x-t} \, dt.$$ It is quite simple to see that these properties immediately imply our goal stated above. Namely, imposing $f=0$ on $w_{-1}^* \setminus \T$ we find $S_nf(x) = S_{k_{-1}}f(x,w_{-1}^*)$ and since all the $j_s$ are pairwise different, we may iterate iv) to obtain $$\big\| S_nf(x) \big\|_\mathrm{X} \, \lesssim \, \big\| S_0f(x,w_L^*) \big\|_\mathrm{X} + M'(f,\varepsilon) \log \log n \, \le \, M'(f,\varepsilon) \log \log n + \lambda.$$ The last estimate follows from $$\|S_0f(x,w_L^*)\|_\mathrm{X} \le \|H_\T^*f(x)\|_\mathrm{X} \le \lambda,$$ since $x \notin \Sigma_{f,\varepsilon}^2$. This proves the desired inequality for $M(f,\varepsilon) = M'(f,\varepsilon) + \lambda$. Recall that $\lambda = \lambda(f,\varepsilon) \to 0^+$ as $f \to 0$ for $\varepsilon$ fixed. Let us then start constructing our family of triplets. As mentioned above, we pick $(k_{-1}, w_{-1}^*) = (n, (-2,2))$. Since $x \in \T \setminus \Sigma_{f,\varepsilon}$, we must have $w_{-1}^* \not\subset \Sigma_{f,\varepsilon}^1$ which gives $$C_{k_{-1}[w_{-1}^*]}^*(f,w_{-1}^*) \, \le \, \max_{\mathrm{grandsons \ of} \ w_{-1}^*} \ \mean_{w'} \|f(x)\|_\mathrm{X} \, dx \, \le \, \lambda.$$ In particular, there must exists $j_{-1} \ge 1$ such that $$2^{-j_{-1}} \lambda \, < \, C_{k_{-1}[w_{-1}^*]}^*(f,w_{-1}^*) \, \le \, 2^{1-j_{-1}} \lambda.$$ This completes the choice of the first triplet $(j_{-1}, k_{-1}, w_{-1}^*)$. Our construction also permits to form the Carleson partition $\Omega_{2^{1-j_{-1}}\lambda}(k_{-1},w_{-1}^*)$. To construct the next triplet, we first consider all the smoothing intervals which arise from the Carleson partition ---i.e. intervals of the form $(2a-b,b)$ or $(a,2b-a)$ for $(a,b)$ an interval in the partition--- which contain $x$ in their middle half. Note that we can always find at least one such interval. Then, we set $w_{-1}^*(x)$ to be the interval of maximal length among the family of smoothing intervals selected. Now we can define the next triplet. First we take $$w_0^* = w_{-1}^*(x) \quad \mbox{and} \quad k_0 = 4 \frac{k_{-1}[w_0^*]}{|w_0^*|}.$$ Then, $j_0$ is determined by $$2^{-j_0} \lambda < C_{k_0[w_0^*]}^*(f,w_0^*) \le 2^{1-j_0} \lambda$$ since $w_0^* \not\subset \Sigma_{f,\varepsilon}^1$ ensures the existence of such a $j_0 \ge 1$. In general, we may produce $(j_{s+1}, k_{s+1}, w_{s+1}^*)$ from the previous triplet in the exact same manner and we get the formulae
\begin{eqnarray*}
w_{s+1}^* & = & w_s^*(x), \\ [5pt] k_{s+1} & = & 4 \frac{k_s[w_{s+1}^*]}{|w_{s+1}^*|}, \\ j_{s+1} & = & 1 + \Big[ \log_2 \Big( \frac{\lambda}{C_{k_{s+1}[w_{s+1}^*]}^*(f,w_{s+1}^*)} \Big) \Big].
\end{eqnarray*}
The square brackets in the last identity stand for the integer part. Recall again that the argument of the $\log_2$ is greater than or equal to $1$ since $w_{s+1}^* = w_s^*(x) \not\subset \Sigma_{f,\varepsilon}^1$ because $x \notin \Sigma_{f,\varepsilon}$. The process finishes at $s=L$ when $k_L = 0$. Since $\frac14 k_s |w_s^*| \in \N$ we will have $k_s=0$ for the first index $s$ with $|w_s^*| < 4/k_s$. In fact, this must happen sooner or later because $|w_s^*|$ is strictly decreasing and $k_s \le n$ for all $s$. Once we have defined the process, let us prove i), ii), iii) and iv). Our choice of $j_s$ is clearly an integer and to show that $j_{s+1} < j_s$ it suffices to see $$C_{k_{s+1}[w_{s+1}^*]}^*(f,w_{s+1}^*) \, > \, 2^{1 - j_s} \lambda.$$ We always have $|w_{s+1}^*| = |w_s^*(s)| \ge 2/n$. This means that the dyadic son $w'$ of $w_s^*(x)$ belonging to $\Omega_{2^{1-j_s}\lambda}(k_s,w_s^*)$ is not of minimal size. In particular, it must have in turn a dyadic son $w''$ satisfying $$C_{k_s[w'']}(f,w'') \, > \, 2^{1-j_s} \lambda.$$ Recalling that $w''$ is a dyadic grandson of $w_{s+1}^*$ and noting that we have the identity $k_s[w''] = k_s[w_{s+1}^*] = k_{s+1}[w_{s+1}^*] = k_{s+1}[w'']$, we conclude that i) holds. Conditions ii) follow from the construction of the maximal intervals $w^*(x)$ and iii) is trivial. It remains to prove iv), which will be done in two steps
\begin{itemize}
\item[a)] Change of period 
\begin{eqnarray*}
\lefteqn{\hskip-5pt \Big| \big\| S_{k_s}f(x,w_{s}^*) \big\|_{\mathrm{X}} - \big\| S_{k_{s}}f(x,w_{s+1}^*) \big\|_{\mathrm{X}} \Big|} \\ & & \hskip50pt \ = \ O \big( 2^{-j_{s+1}(1 - \frac{1}{1+\delta/2})} M''(f,\varepsilon) \log \log n \big).
\end{eqnarray*}

\vskip3pt

\item[b)] Change of frequency 
$$\Big| \big\| S_{k_s}f(x,w_{s+1}^*) \big\|_{\mathrm{X}} - \big\| S_{k_{s+1}}f(x,w_{s+1}^*) \big\|_{\mathrm{X}} \Big| \ = \ O \big( 2^{-j_{s+1}} \lambda \big).$$
\end{itemize}
For the change of frequency, we note that
\begin{eqnarray*}
\lefteqn{\hskip-1pt \Big| \big\| S_{k_s}f(x,w_{s+1}^*) \big\|_{\mathrm{X}} - \big\| S_{k_{s+1}}f(x,w_{s+1}^*) \big\|_{\mathrm{X}} \Big|} \\ & \le & \Big\| e^{2\pi i k_s x} S_{k_s}f(x,w_{s+1}^*) - e^{2\pi i k_{s+1} x} S_{k_{s+1}}f(x,w_{s+1}^*) \Big\|_{\mathrm{X}} \\ & = & \Big\| \int_{w_{s+1}^*} \frac{f(t) (e^{2\pi i k_s(x-t)} - e^{2 \pi i k_{s+1} (x-t)})}{x-t} \, dt \Big\|_{\mathrm{X}} \\ & \le & \sum_{\mathrm{grandsons}} \hskip1pt \Big\| \int_{w'} \phi_s(x-t) f(t) e^{2 \pi i k_{s+1} (x-t)} \, dt \Big\|_{\mathrm{X}},
\end{eqnarray*}
with $\phi_s(x) = \frac{1}{x} (e^{2\pi i (k_s - k_{s+1})x} - 1)$. Arguing as in the proof of Lemma \ref{Lemmagoodpart} $$\phi_s(x) = \sum_{\mu \in \Z} \gamma_\mu \exp(- 2 \pi i \mu x /3|w'|)$$ with $$(1+\mu^2)|\gamma_\mu| \lesssim \|\phi_s\|_\infty + |w'|^2 \|\phi_s''\|_\infty \lesssim \frac{1}{|w'|}.$$ The last inequality follows from $|w'| |k_s - k_{s+1}| \lesssim 1$. Implementing this above yields 
\begin{eqnarray*}
\lefteqn{\hskip-9pt \Big| \big\| S_{k_s}f(x,w_{s+1}^*) \big\|_{\mathrm{X}} - \big\| S_{k_{s+1}}f(x,w_{s+1}^*) \big\|_{\mathrm{X}} \Big|} \\ & \lesssim & \sum_{\mathrm{grandsons}} C_{k_{s+1}[w']}(f,w') \ \lesssim \ C_{k_{s+1}[w_{s+1}^*]}^*(f,w_{s+1}^*).
\end{eqnarray*}
We know the term on the right is bounded above by $2^{1 - j_{s+1}} \lambda$. Therefore, the proof of b) is complete. 
For the change of interval, we use that $x \notin \Sigma_{f,\varepsilon}^3 \cup \Sigma_{f,\varepsilon}^4$ as follows
\begin{eqnarray*}
\lefteqn{\hskip-5pt \Big| \big\| S_{k_s}f(x,w_{s}^*) \big\|_{\mathrm{X}} - \big\| S_{k_{s}}f(x,w_{s+1}^*) \big\|_{\mathrm{X}} \Big|} \\ & \le & \Big\| \int_{w_s^* \setminus w_s^*(x)} \frac{f(t) e^{- 2 \pi i k_s t}}{x-t} \, dt \Big\|_\mathrm{X} \\ & \le & \Big\| \int_{w_s^* \setminus w_s^*(x)} \frac{g_{k_s, 2^{1-j_s}\lambda}(t)}{x-t} \, dt \Big\|_\mathrm{X} + \Big\| \int_{w_s^* \setminus w_s^*(x)} \frac{b_{k_s, 2^{1-j_s}\lambda}(t)}{x-t} \, dt \Big\|_\mathrm{X} \ = \ \mathrm{A} + \mathrm{B}.
\end{eqnarray*}
Since we clearly have $$(k_s, w_s^*) \in \bigcup_{m \ge \frac{1}{\varepsilon} C(\lambda,\delta)} \mathcal{A}_{e^m} \quad \mbox{for every} \quad -1 \le s \le L$$ and $x \notin \Sigma_{f,\varepsilon}^3$, we conclude $$\mathrm{A} \, \le \, 2 \big\| H_{w_s^*}^*g_{k_s,2^{1-j_s}\lambda}(x) \big\|_\mathrm{X} \, \le \, R(\lambda) (2^{1-j_s} \lambda)^{1- \frac{1}{1+\delta/2}} \log\log n.$$ On the other hand, we claim that $$\mathrm{B} \, \lesssim \, 2^{-j_s} \lambda \, \Delta_{\Omega_{2^{1-j_s} \lambda}(k_s,w_s^*)}(x).$$ Assuming the claim, we may argue as for the $\mathrm{A}$-term and obtain the exact same upper bound for $\mathrm{B}$. The only difference is that we now have to use the fact that $x \notin \Sigma_{f,\varepsilon}^4$. Since $j_s > j_{s+1}$, we find $$M''(f,\varepsilon) = \lambda^{1 - \frac{1}{1+\delta/2}} R(\lambda) \to 0^+ \quad \mbox{as} \quad f \to 0 \quad \mbox{for fixed $\varepsilon$},$$ since $\lambda = \frac{1}{\varepsilon} \|f\|_{L (\log L)^{1+\delta}}$. Moreover, according to b) we get $M'(f,\varepsilon) = M''(f,\varepsilon) + \lambda$ which again goes to $0$ with $f$ for $\varepsilon$ fixed. Consequently, it just remains to prove our claim for the term $\mathrm{B}$. A moment of thought gives ---recalling the definition of the interval $w^*(x)$ associated to a Carleson partition $\Omega_\lambda(k,w^*)$--- the following two properties associated to the maximality of $w^*(x)$
\begin{itemize}
\item $w^*(x)$ is a union of intervals in $\Omega_\lambda(k,w^*)$,

\item Given $w' \in \Omega_\lambda(k,w^*) \setminus w^*(x)$ we have $\mathrm{dist}(x,w') \ge \frac12 |w'|$.
\end{itemize}
Using the mean zero of $b_{k_s, 2^{1-j_s} \lambda}$ on the intervals of $\Omega_{2^{1-j_s} \lambda}(k_s,w_s^*)$, we find 
\begin{eqnarray*}
\mathrm{B} & \le & \sum_{w' \in \Omega_{2^{1-j_s} \lambda}(k_s,w_s^*) \setminus w_s^*(x)} \Big\| \int_{w'} \frac{t-\mathrm{c}(w')}{(x-t)(x-\mathrm{c}(w'))} f(t) e^{-2 \pi i k_st} \, dt \Big\|_\mathrm{X} \\ & + & \sum_{w' \in \Omega_{2^{1-j_s} \lambda}(k_s,w_s^*) \setminus w_s^*(x)} \Big\| \int_{w'} \frac{t-\mathrm{c}(w')}{(x-t)(x-\mathrm{c}(w'))} g_{k_s, 2^{1-j_s} \lambda} \, dt \Big\|_\mathrm{X} \ = \ \mathrm{B}_1 + \mathrm{B}_2.
\end{eqnarray*}
According to Lemma \ref{Lemmagoodpart} and the fact that $\mathrm{dist}(x,w') \ge \frac12 |w'|$ we immediately get $\mathrm{B}_2 \lesssim 2^{-j_s} \lambda \Delta_{\Omega_{2^{1-j_s} \lambda}(k_s,w_s^*)}(x)$. Using one more time the Fourier expansion argument in the proof of Lemma \ref{Lemmagoodpart} for the function $$\phi_{w'}(t) = \frac{t - \mathrm{c}(w')}{(x-t)(x-\mathrm{c}(w'))} \exp \Big( - 2 \pi i \big(k_s - \frac{k_s[w']}{|w'|} \big) t \Big)$$ and the Whitney property $\mathrm{dist}(x,w') \ge \frac12 |w'|$, we easily get $$\mathrm{B}_1 \lesssim \sum_{w' \in \Omega_{2^{1-j_s} \lambda}(k_s,w_s^*)} \frac{|w'|^2 C_{k_s[w']}(f,w')}{(x - \mathrm{c}(w'))^2 + |w'|^2} \, \le \, 2^{1-j_s} \lambda \Delta_{\Omega_{2^{1-j_s} \lambda}(k_s,w_s^*)}(x).$$ This completes the argument for the claim, and Theorem A is proved. \hspace*{\fill} $\blacksquare$ \vskip0.2cm
 
\vskip10pt

\noindent \textbf{Acknowledgement.} The authors of this paper were supported in part by the ERC Grant StG-256997-CZOSQP (European Union), the MEC Grant MTM-2010-16518 (Spain) and ANR (France).

\bibliographystyle{amsplain}

\vskip20pt

\hfill \noindent \textbf{Javier Parcet} \\
\null \hfill Instituto de Ciencias Matem{\'a}ticas \\ \null \hfill
CSIC-UAM-UC3M-UCM \\ \null \hfill Consejo Superior de
Investigaciones Cient{\'\i}ficas \\ \null \hfill C/ Nicol\'as Cabrera 13-15.
28049, Madrid. Spain \\ \null \hfill\texttt{javier.parcet@icmat.es}

\

\hfill \noindent \textbf{Fernando Soria} \\
\null \hfill Instituto de Ciencias Matem{\'a}ticas \\ \null \hfill
CSIC-UAM-UC3M-UCM \\ \null \hfill Universidad Aut\'onoma de Madrid 
\\ \null \hfill C/ Nicol\'as Cabrera 13-15.
28049, Madrid. Spain \\ \null \hfill\texttt{fernando.soria@uam.es}

\

\hfill \noindent \textbf{Quanhua Xu} \\
\null \hfill School of Mathematics and Statistics 
\\ \null \hfill Wuhan University
\\ \null \hfill Wuhan, 430072, China
\\ \null \hfill Laboratoire de Math{\'e}matiques 
\\ \null \hfill Universit{\'e} de France-Comt{\'e}
\\ \null \hfill 16 Route de Gray, 25030 Besan\c{c}on Cedex, France 
\\ \null \hfill\texttt{qxu@univ-fcomte.fr}

\end{document}